\documentclass{amsart}
\usepackage{amstext}
\usepackage{amsthm}
\usepackage{amssymb}

\usepackage{hyperref}

\makeatletter
\numberwithin{equation}{section}
\numberwithin{figure}{section}
\theoremstyle{plain}
\newtheorem{thm}{Theorem}

\theoremstyle{remark}
\newtheorem*{rmk}{Remark}

\allowdisplaybreaks
\usepackage[sort&compress,square,numbers]{natbib}
\makeatother

\usepackage{tikz}
\usetikzlibrary{matrix,positioning}

\begin{document}

\title[A new approach to Catalan numbers using differential equations]{A new approach to Catalan numbers using differential equations}

\author{Taekyun Kim}
\address{Department of Mathematics, Kwangwoon University, Seoul 139-701, Republic
	of Korea}
\email{tkkim@kw.ac.kr}

\author{Dae San Kim}
\address{Department of Mathematics, Sogang University, Seoul 121-742, Republic
	of Korea}
\email{dskim@sogang.ac.kr}
\begin{abstract}
In this paper, we introduce two differential equations arising from the generating function of the Catalan numbers which are `inverses' to each other in some sense. From these differential equations, we obtain some new and explicit identities for Catalan and higher-order Catalan numbers. In addition, by other means than differential equations we also derive some interesting identities involving Catalan numbers which are of arithmetic and combinatorial nature.   
\end{abstract}

 \medskip

\maketitle
\section{Introduction}
The Catalan numbers $C_n$ were first introduced by the Mongolian mathematician Ming Antu in around 1730, even though they were named after the Belgian mathematician Eug{\`e}ne Charles Catalan (1814-1894). Indeed, Ming Antu obtained a number of trigonometric expressions involving Catalan numbers such as
\begin{equation*}\begin{split}
\sin 2\theta = 2 \sin \theta - \sum_{n=1}^\infty \frac{C_{n-1}}{4^{n-1}}\sin^{2n+1}\theta
=2 \sin \theta - \sin^3 \theta - \tfrac{1}{4} \sin^5 \theta - \tfrac{1}{8} \sin^7 \theta-\cdots, \,\,\, 
\end{split}\end{equation*}
(\textnormal{see} \cite{key-2,key-3,key-4,key-6,key-10,key-11,key-12,key-13,key-14}).
The Catalan numbers  can be given explicitly in terms of binomial coefficients. Namely, for $n \geq 0$,
\begin{equation}\begin{split}\label{01}
C_n = \frac{1}{n+1} {2n \choose n} = \prod_{k=2}^n \frac{n+k}{k}.
\end{split}\end{equation}
They satisfy the recurrence relations
\begin{equation}\begin{split}\label{02}
C_0=1, C_n = \sum_{m=0}^{n-1}C_m C_{n-1-m}, \quad (n \geq 1).
\end{split}\end{equation}
The Catalan numbers are also given by the generating function
\begin{equation}\begin{split}\label{03}
\frac{2}{1+\sqrt{1-4t}} = \sum_{n=0}^\infty C_n t^n = \sum_{n=0}^\infty \frac{1}{n+1} {2n \choose n} t^n.
\end{split}\end{equation}
The Catalan numbers form the sequence of positive integers
\begin{equation*}\begin{split}
1,\,\, 1, \,\,2,\,\, 5,\,\,  14,\,\, 42,\,\, 132,\,\, 429,\,\, 1430,\,\, 4862,\,\, 16796,\,\, 58786, \,\,208012, \,\,\cdots
\end{split}\end{equation*}
which is asymptotic to $\dfrac{4^n}{n^{\tfrac{3}{2}} \sqrt{\pi}}$, as $n$ tends to $\infty$, and appears in various counting problems. For example, $C_n$ is the number of Dyck words of length $2n$, the number of balanced $n$ pairs of parentheses, the number of mountain ranges you can form with $n$ upstrokes and  downstrokes that all stay above the original line, the number of diagonal-avoiding paths of length $2n$ from the upper left corner to the lower right corner in a grid of $n \times n$ squares, and the number of ways $n+1$ factors can be completely parenthesized (see \cite{key-1,key-2,key-3,key-4,key-10,key-11}). 

 It is also the number of ways an $(n+2)-$gon can be cut into $n$ triangles, the number of permutations of $\{ 1,2,\cdots,n \}$ that avoid the pattern 123, the number of ways to tile a stairstep shape of height $n$ with $n$ rectangles, etc (see \cite{key-10,key-11,key-12,key-13}).
  
  In \cite{key-8}, T. Kim initiated a fascinating idea of using ordinary differential equations 
as a method of obtaining new identities for special polynomials and numbers. Namely, a family of nonlinear differential equations  were derived, which are indexed by positive integers and satisfied by the generating function of the Frobenius-Euler numbers. Then, they were used in order to obtain an interesting identity, expressing higher-order Frobenius-Euler numbers in terms of (ordinary) Frobenius-Euler numbers (see \cite{key-5,key-7,key-9,key-15}). 

This method turned out to be very fruitful and can be applied to many interesting special polynomials and numbers  (see \cite{key-7,key-8,key-9}). For example, linear differential equations are derived for Bessel polynomials, Changhee polynomials, actuarial polynomials, Meixner polynomials of the first kind, Poisson-Charlier polynomials, Laguerre polynomials, Hermite polynomials,
and Stirling polynomials, while nonlinear ones are obtained for Bernoulli numbers of the second, Boole numbers, Chebyshev polynomials of the first, second, third, and fourth kind, degenerate Euler numbers, degenerate Eulerian polynomials, Korobov numbers, and Legendre polynomials (see \cite{key-1,key-5,key-7,key-9,key-15}).

  In this paper, we introduce two differential equations arising from the generating function of the Catalan numbers which are `inverses' to each other in some sense. From these differential equations, we obtain some new and explicit identities for Catalan and higher-order Catalan numbers. In addition, by other means than differential equations we also derive some interesting identities involving Catalan numbers which are of arithmetic and combinatorial nature.   

\section{Differential equations associated with Catalan numbers}
Let
\begin{equation}\begin{split}\label{04}
C=C(t) = \frac{2}{1+\sqrt{1-4t}}
\end{split}\end{equation}
Then, by \eqref{04}, we have
\begin{equation}\begin{split}\label{05}
C^{(1)} =& \dfrac{d}{dt} C(t) = (1-4t)^{-\tfrac{1}{2}} 4(1+ \sqrt{1-4t})^{-2}\\
=&(1-4t)^{-\tfrac{1}{2}} C^2,
\end{split}\end{equation}
and
\begin{equation}\begin{split}\label{06}
C^{(2)} =& \dfrac{d}{dt}C^{(1)} = 2 (1-4t)^{-\tfrac{3}{2}} C^2 + 2(1-4t)^{-\tfrac{1}{2}} C C^{(1)} \\
=& 2 (1-4t)^{-\tfrac{3}{2}} C^2 + 2 (1-4t)^{-\tfrac{1}{2}} C \{(1-4t)^{-\tfrac{1}{2}} C^2 \}\\
=& 2(1-4t)^{-\tfrac{3}{2}} C^2 + 2(1-4t)^{-\tfrac{2}{2}} C^3
\end{split}\end{equation}
So, we are led to put
\begin{equation}\begin{split}\label{07}
C^{(N)} = \sum_{i=1}^N a_i(N) (1-4t)^{-\tfrac{2N-i}{2}} C^{i+1},
\end{split}\end{equation}
where $N=1,2,3,\cdots$.
From \eqref{07}, we obtain
\begin{equation}\begin{split}\label{08}
C^{(N+1)} =& \dfrac{d}{dt}C^{(N)}\\
=& \sum_{i=1}^N 2(2N-i)a_i(N)(1-4t)^{-\tfrac{2N+2-i}{2}} C^{i+1}\\
&+\sum_{i=1}^N (i+1) a_i(N) (1-4t)^{-\tfrac{2N-i}{2}}C^i C^{(1)}\\
=& \sum_{i=1}^N 2(2N-i) a_i(N) (1-4t)^{-\tfrac{2N+2-i}{2}}  C^{i+1}\\
&+ \sum_{i=1}^N (i+1)a_i(N) (1-4t)^{-\tfrac{2N-i}{2}}C^i \{ (1-4t)^{-\tfrac{1}{2}} C^2 \}\\
=& \sum_{i=1}^N 2(2N-i) a_i(N) (1-4t)^{-\tfrac{2N+2-i}{2}}  C^{i+1}\\
&+ \sum_{i=1}^N (i+1)a_i(N) (1-4t)^{-\tfrac{2N+1-i}{2}} C^{i+2}\\
=& \sum_{i=1}^N 2(2N-i) a_i(N) (1-4t)^{-\tfrac{2N+2-i}{2}}  C^{i+1}\\
&+ \sum_{i=2}^{N+1} i a_{i-1}(N) (1-4t)^{-\tfrac{2N+2-i}{2}}C^{i+1}.
\end{split}\end{equation}
On the other hand, replacing $N$ by $N+1$ in \eqref{07}, we get
\begin{equation}\begin{split}\label{09}
C^{(N+1)} = \sum_{i=1}^{N+1} a_i(N+1)(1-4t)^{-\tfrac{2N+2-i}{2}}C^{i+1}.
\end{split}\end{equation}
From \eqref{08} and \eqref{09}, we can derive the following recurrence relations:
\begin{align}
\label{10}&a_1(N+1)=2(2N-1)a_1(N),\\
\label{11}&a_{N+1}(N+1)=(N+1)a_N(N),
\end{align}
and
\begin{equation}\begin{split}\label{12}
a_i(N+1)=ia_{i-1}(N) + 2(2N-i) a_i(N), \,\, (2 \leq i \leq N).
\end{split}\end{equation}
In addition, from \eqref{05} and \eqref{07}, we observe that
\begin{equation}\begin{split}\label{13}
a_1(1)(1-4t)^{-\tfrac{1}{2}} C^2 = C^{(1)} = (1-4t)^{-\tfrac{1}{2}} C^2.
\end{split}\end{equation}
Thus, by \eqref{13}, we get
\begin{equation}\begin{split}\label{14}
a_1(1) = 1.
\end{split}\end{equation}
In Below, for any positive integer $N$, $(2N-1)!!$ will denote
\begin{equation}\begin{split}\label{15}
(2N-1)!! = (2N-1)(2N-3) \cdots 1.
\end{split}\end{equation}
From \eqref{10} and \eqref{14}, we note that
\begin{equation}\begin{split}\label{16}
a_1(N+1)=& 2(2N-1)a_1(N) = 2^2 (2N-1)(2N-3)a_1(N-1)\\
=& \cdots \\
=& 2^N(2N-1)(2N-3)\cdots 1a_1(1)\\
=& 2^N(2N-1)!!
\end{split}\end{equation}
and
\begin{equation}\begin{split}\label{17}
a_{N+1}(N+1)=& (N+1)a_N(N) = (N+1)Na_{N-1}(N-1) \\
=& \cdots \\
=& (N+1)N \cdots 2 a_1(1) = (N+1)!
\end{split}\end{equation}
In the following, we will use the notations:
\begin{equation}\begin{split}\label{18}
(x;\alpha)_n = x(x-\alpha ) \cdots ( x- (n-1)\alpha ),\quad \textnormal{for}\,\, n \geq 1,
\end{split}\end{equation}
and
\begin{equation*}\begin{split}
(x; \alpha)_0 = 1.
\end{split}\end{equation*}
For $i=2$ in \eqref{12}, we have
\begin{equation}\begin{split}\label{19}
a_2(N+1) =& 2a_1(N)+2(2N-2)a_2(N) \\
=& 2a_1(N)+2(2N-2)\big(2a_1(N-1)+2(2N-4)a_2(N-1)\big)\\
=& 2\big( a_1(N) + 2(2N-2) a_1(N-1) \big) + 2^2 (2N-2) (2N-4) a_2(N-1)\\
=& 2\big( a_1(N) + 2(2N-2) a_1(N-1) \big) \\
&+ 2^2 (2N-2) (2N-4)\big( 2a_1(N-2)+2(2N-6)a_2(N-2)\big)\\
=& 2\big( a_1(N) + 2(2N-2) a_1(N-1) \big) + 2^2 (2N-2) (2N-4)a_1(N-2)\big)\\
&+2^3 (2N-2)(2N-4)(2N-6)a_2(N-2)\\
=&\,\,\cdots\\
=&2\sum_{k=0}^{N-2} 2^k (2N-2;2)_k a_1(N-k) + 2^{N-1}(2N-2;2)_{N-1}a_2(2)\\
=&2\sum_{k=0}^{N-1}2^k (2N-2;2)_k a_1(N-k).
\end{split}\end{equation}
Proceeding analogously to the case of $i=2$, for $i=3$ and 4, we obtain
\begin{align}
\label{20}&a_3(N+1)=3\sum_{k=0}^{N-2}2^k(2N-3;2)_ka_2(N-k),\\
\label{21}&a_4(N+1)=4\sum_{k=0}^{N-3}2^k(2N-4;2)_ka_3(N-k).
\end{align}
Continuing this process, we can deduce that
\begin{equation}\begin{split}\label{22}
a_i(N+1)=i \sum_{k=0}^{N-i+1}2^k (2N-i;2)_k a_{i-1}(N-k), \quad \textnormal{for} \,\, 2 \leq i \leq N.
\end{split}\end{equation}
Now, we give explicit expressions for $a_i(N+1)\,\,(2 \leq i \leq N).$ From \eqref{16} and \eqref{19}, we have
\begin{equation}\begin{split}\label{23}
a_2(N+1) =& 2 \sum_{k_1=0}^{N-1} 2^{k_1} (2N-2;2)_{k_1} a_1(N-k_1) \\
=& 2 \sum_{k_1=0}^{N-1} 2^{k_1}(2N-2;2)_{k_1} 2^{N-k_1-1}(2N-2k_1-3)!!\\
=& 2! 2^{N-1} \sum_{k_1=0}^{N-1} (2N-2;2)_{k_1} (2N-2k_1 -3)!!.
\end{split}\end{equation}
Also, from \eqref{20} and \eqref{23}, we get
\begin{equation}\begin{split}\label{24}
a_3(N+1)=&3\sum_{k_2=0}^{N-2}2^{k_2}(2N-3;2)_{k_2} a_2(N-k_2) \\
=& 3 \sum_{k_2=0}^{N-2} 2^{k_2} (2N-3;2)_{k_2} 2^{N-k_2-1} \\&\times \sum_{k_1=0}^{N-2-k_2} (2N-2k_2-4;2)_{k_1} (2N-2k_1-2k_1-5)!! \\
=& 3! 2^{N-2} \sum_{k_2=0}^{N-2}\sum_{k_1=0}^{N-2-k_2}(2N-3;2)_{k_2} (2N-2k_2-4;2)_{k_1}\\ 
&\times (2N-2k_1-2k_1-5)!!.
\end{split}\end{equation}
Continuing this process, we can deduce that
\begin{equation}\begin{split}\label{25}
&a_i(N+1)= 2^{N-i+1} i! \sum_{k_{i-1}=0}^{N-i+1} \sum_{k_{i-2}=0}^{N-i+1-k_{i-1}} \cdots \sum_{k_1=0}^{N-i+1-k_{i-1}-\cdots-k_2} (2N-i;2)_{k_{i-1}}\\
&\,\,\,\,(2N-2k_{i-1}-i-1;2)_{k_{i-2}} \times \cdots \times (2N-2k_{i-1}-\cdots-2k_2-2i+2;2)_{k_1}\\
&\,\,\,\,\times (2N-2k_{i-1}-\cdots-2k_1 -2i +1)!!\\
&= 2^{N-i+1}i! \sum_{k_{i-1}=0}^{N-i+1} \sum_{k_{i-2}=0}^{N-i+1-k_{i-1}} \cdots \sum_{k_1=0}^{N-i+1-k_{i-1}-\cdots-k_2} \\
&\quad \prod_{l=1}^{i-1} (2N-2 \sum_{j=l+1}^{i-1} k_j -2i+1+l;2)_{k_l} \times (2N-2 \sum_{j=1}^{i-1} k_j -2i+1)!!,
\end{split}\end{equation}
for $2 \leq i \leq N$.

\begin{rmk}
We note here that \eqref{25} is also valid for $i=N+1$. 
\end{rmk}

Therefore, from \eqref{16} and \eqref{25}, we obtain the following theorem.
\begin{thm}\label{thm:1}
The family of differential equations 
\begin{equation}\begin{split}\label{26}
C^{(N)} = \sum_{i=1}^N a_i(N) (1-4t)^{-\tfrac{2N-i}{2}} C^{i+1} \quad (N=1,2,3,\cdots)
\end{split}
\end{equation}
have a solution 
\begin{equation}
C=C(t)= \tfrac{2}{1+\sqrt{1-4t}}\nonumber, 
\end{equation}
where
\begin{align*}
a_1(N)& =2^{N-1}(2N-3)!!, \\
a_i(N)& =2^{N-i}i! \sum_{k_{i-1}=0}^{N-i} \sum_{k_{i-2}=0}^{N-i-k_{i-1}} \cdots \sum_{k_1=0}^{N-i-k_{i-1}-\cdots-k_2} \qquad\qquad\qquad  \\
&\times \prod_{l=1}^{i-1} (2N-2 \sum_{j=l+1}^{i-1} k_j -2i-1+l;2)_{k_l} \\ 
&\times (2N-2 \sum_{j=1}^{i-1} k_j -2i-1)!!.
\end{align*}
\end{thm}

We recall that the Catalan numbers $C_n$ are defined by the generating funcion
\begin{equation}\begin{split}\label{27}
C=C(t) = \frac{2}{1+\sqrt{1-4t}} = \sum_{n=0}^\infty C_n t^n.
\end{split}\end{equation}
More generally, the higher-order Catalan numbers $C_n^{(r)}$ of order $r$ are given by
\begin{equation}\begin{split}\label{28}
\left( \frac{2}{1+\sqrt{1-4t}} \right)^r = \sum_{n=0}^\infty C_n^{(r)} t^n.
\end{split}\end{equation}
On the one hand, from \eqref{27}, we have
\begin{equation}\begin{split}\label{29}
C^{(N)} =& \sum_{n=N}^\infty C_n(n)_N t^{n-N}\\
=& \sum_{n=0}^\infty C_{n+N} (n+N)_N t^n,
\end{split}\end{equation}
where 
\begin{equation}\begin{split}\label{30}
(x)_n = x(x-1) \cdots (x-n+1), \quad \textnormal{for}\,\, n \geq 1, \,\, (x)_0=1.
\end{split}\end{equation}
On the other hand, by Theorem \ref{thm:1}, we have
\begin{equation}\begin{split}\label{31}
C^{(N)} =& \sum_{i=1}^N a_i(N) (1-4t)^{-\tfrac{2N-i}{2}} C^{i+1}\\
=& \sum_{i=1}^N a_i(N) \sum_{m=0}^\infty { \tfrac{2N-i}{2}+m-1 \choose m } 4^m t^m \sum_{l=0}^\infty C_l^{(i+1)} t^l \\
=& \sum_{i=1}^N a_i(N) \sum_{n=0}^\infty \sum_{m=0}^n 4^m  { \tfrac{2N-i}{2}+m-1 \choose m } C_{n-m}^{(i+1)} t^n\\
=& \sum_{n=0}^\infty \left( \sum_{i=1}^N \sum_{m=0}^n 4^m  { \tfrac{2N-i}{2}+m-1 \choose m }  a_i(N) C_{n-m}^{(i+1)} \right) t^n.
\end{split}\end{equation}
Comparing \eqref{29} with \eqref{31}, we get the following Theorem.

\begin{thm}\label{thm:2}
For $n=0,1,2,\cdots,$ and $N=1,2,3,\cdots,$
\begin{equation*}\begin{split}
C_{n+N}=& \frac{1}{(n+N)_N} \sum_{i=1}^N \sum_{m=0}^n 4^m  { \tfrac{2N-i}{2}+m-1 \choose m }  a_i(N) C_{n-m}^{(i+1)}
\end{split}\end{equation*}
where $a_i(N)$'s are as in Theorem \ref{thm:1}.
\end{thm}

\section{Inverse differential equations associated with Catalan numbers}
Here we will derive ``inverse'' differential equations to the ones obtained in Section 2. With $C=C(t)$ as in \eqref{04}, we have
\begin{equation}\begin{split}\label{32}
C^{(1)}=(1-4t)^{-\tfrac{1}{2}} C^2,
\end{split}\end{equation}
and
\begin{equation}\begin{split}\label{33}
C^2 = (1-4t)^{\tfrac{1}{2}}C^{(1)}.
\end{split}\end{equation}
Differentiating both sides of \eqref{33}, we get
\begin{equation}\begin{split}\label{34}
2CC^{(1)} = -2(1-4t)^{-\tfrac{1}{2}}C^{(1)} + (1-4t)^{\tfrac{1}{2}}C^{(2)}.
\end{split}\end{equation}
Substituting \eqref{32}, into \eqref{34}, we obtain
\begin{equation}\begin{split}\label{35}
2C^3 = -2C^{(1)} + (1-4t)C^{(2)}.
\end{split}\end{equation}
Differentiating both sides of \eqref{35}, we have
\begin{equation}\begin{split}\label{36}
3! C^2 C^{(1)} = -6C^{(2)} + (1-4t)C^{(3)}
\end{split}\end{equation}
Substituting \eqref{32} into \eqref{36}, we get
\begin{equation}\begin{split}\label{37}
3! C^4 = -6(1-4t)^{\tfrac{1}{2}}C^{(2)} +(1-4t)^{\tfrac{3}{2}} C^{(3)}.
\end{split}\end{equation}
So we are led to put
\begin{equation}\begin{split}\label{38}
N! C^{N+1} = \sum_{i=0}^{\left[ \frac{N}{2}\right]}b_i(N) (1-4t)^{\tfrac{N}{2}-i} C^{(N-i)} \quad (N=1,2,3,\cdots).
\end{split}\end{equation}
Here $[x]$ denote the greatest integer not exceeding $x$. Differentiation of both sides of \eqref{38}, gives
\begin{equation}\begin{split}\label{39}
(N+1)! C^N C^{(1)} =& \sum_{i=0}^{\left[ \frac{N}{2}\right]} -4(\tfrac{N}{2}-i) b_i(N) (1-4t)^{\tfrac{N}{2}-i-1} C^{(N-i)}\\
&+  \sum_{i=0}^{\left[ \frac{N}{2}\right]}b_i(N) (1-4t)^{\tfrac{N}{2}-i} C^{(N+1-i)}\\
=& \sum_{i=1}^{\left[ \frac{N}{2}\right]+1} -4(\tfrac{N}{2}+1-i) b_{i-1}(N) (1-4t)^{\tfrac{N}{2}-i} C^{(N+1-i)}\\
&+  \sum_{i=0}^{\left[ \frac{N}{2}\right]}b_i(N) (1-4t)^{\tfrac{N}{2}-i} C^{(N+1-i)}.
\end{split}\end{equation}
Substituting \eqref{32} into \eqref{39}, we obtain
\begin{equation}\begin{split}\label{40}
(N+1)! C^{N+1} =& \sum_{i=1}^{\left[ \frac{N}{2}\right]+1} -4( \tfrac{N}{2} +1 -i) b_{i-1}(N)(1-4t)^{\frac{N+1}{2}-i} C^{(N+1-i)}\\
&+ \sum_{i=0}^{ \left[ \frac{N}{2} \right]} b_i(N) (1-4t)^{\tfrac{N+1}{2}-i} C^{(n+1-i)}.
\end{split}\end{equation}
Also, by replacing $N$ by $N+1$ in \eqref{38}, we get
\begin{equation}\begin{split}\label{41}
(N+1)! C^{N+2} = \sum_{i=0}^{\left[ \frac{N+1}{2} \right]} b_i(N+1) (1-4t)^{\tfrac{N+1}{2}-i} C^{(N+1-i)}.
\end{split}\end{equation}
Comparing \eqref{40} with \eqref{41}, we have the following recurrence relations. Here we need to consider the even and odd cases of $N$ separately. The details are left to the reader.
\begin{align}
\label{42} &b_0(N+1)=b_0(N),\\
\label{43} &b_i(N+1)=-4( \tfrac{N}{2}+1-i) b_{i-1}(N) + b_i(N),\quad \textnormal{for} \,\,\, 1 \leq i \leq \left[ \tfrac{N+1}{2} \right].
\end{align}
From \eqref{33} and \eqref{38}, we have 
\begin{equation}\begin{split}\label{44}
C^2 = b_0(1) (1-4t)^{\tfrac{1}{2}}C^{(1)} = (1-4t)^{\tfrac{1}{2}}C^{(1)}.
\end{split}\end{equation}
Thus, from \eqref{44}, we get
\begin{equation}\begin{split}\label{45}
b_0(1) = 1.
\end{split}\end{equation}
From \eqref{42}, we easily obtain
\begin{equation}\begin{split}\label{46}
b_0(N+1)=b_0(N)=\cdots=b_0(1) =1.
\end{split}\end{equation}
The equation in \eqref{43} can be rewritten as
\begin{equation}\begin{split}\label{47}
b_i(N+1) = -2(N+2-2i) b_{i-1}(N) + b_i(N)
\end{split}\end{equation}
To proceed further, we define
\begin{align}
\label{48} &S_{N,1} = N+(N-1)+\cdots+1,\\
\label{49} &S_{N,j} = NS_{N+1,j-1} + (N-1)S_{N,j-1} + \cdots + 1 S_{2,j-1} \quad (j \geq 2).
\end{align}
Now,
\begin{equation}\begin{split}\label{50}
b_1(N+1) =& -2Nb_0(N)+b_1(N)
\\=& -2N+b_1(N)
\\=& -2N-2(N-1)+b_1(N-1)
\\=& \cdots
\\=& -2(N+(N-1)\cdots+1) + b_1(1)
\\=& -2 S_{N,1},
\end{split}\end{equation}
and
\begin{equation}\begin{split}\label{51}
b_2(N+1)=&-2(N-2)b_1(N) +b_2(N)
\\=&(-2)^2 (N-2) S_{N-1,1} + b_2(N)
\\=&(-2)^2 \big( (N-2)S_{N-1,1} + (N-3)S_{N-2,1} \big) + b_2(N-1)
\\=&\cdots
\\=&(-2)^2 \big( (N-2)S_{N-1,1} + (N-3)S_{N-2,1} + \cdots + 1S_{2,1} \big) + b_2(3)
\\=&(-2)^2 S_{N-2,2}.
\end{split}\end{equation}
Similarly to the cases of $i=1$ and 2, for $i=3$, we get
\begin{equation}\begin{split}\label{52}
b_3(N+1) = (-2)^3 S_{N-4,3}.
\end{split}\end{equation}
Thus we can deduce that, for $1 \leq i \leq \left[ \tfrac{N+1}{2} \right],$
\begin{equation}\begin{split}\label{53}
b_i(N+1) = (-2)^i S_{N+2-2i, i}.
\end{split}\end{equation}
Here, from \eqref{46} and \eqref{53}, we obtain the following Theorem.

\begin{thm}\label{thm:3}
The following family of differential equations
\begin{equation}\begin{split}\label{54}
N! C^{N+1} = \sum_{i=0}^{ \left[ \tfrac{N}{2} \right]} b_i(N)(1-4t)^{\tfrac{N}{2}-i} C^{(N-i)} \quad (N=1,2,3\cdots)
\end{split}\end{equation}
have a solution 
\begin{equation}
C=C(t)=\frac{2}{1+\sqrt{1-4t}},\nonumber 
\end{equation}
where 
\begin{equation}
b_0(N)=1, b_i(N) = (-2)^i S_{N+1-2i,i}, \quad ( 1 \leq i \leq \left[ \tfrac{N}{2} \right])\nonumber.
\end{equation}
\end{thm}

Now, we would like to give an application of the result in Theorem \ref{thm:3}. From \eqref{54}, we have the following
\begin{equation}\begin{split}\label{55}
N! \sum_{k=0}^\infty C_k^{(N+1)} t^k=& \sum_{i=0}^{ \left[ \tfrac{N}{2} \right]}  b_i(N) \sum_{l=0}^\infty { \tfrac{N}{2}-i \choose l} (-4t)^l\\
&\times \sum_{m=0}^\infty C_{m+N-i}(m+N-i)_{N-i} t^m \\
=& \sum_{i=0}^{ \left[ \tfrac{N}{2} \right]} b_i(N) \sum_{k=0}^\infty \sum_{m=0}^k {\tfrac{N}{2}-i \choose k-m} (-4)^{k-m} \\
&\times C_{m+N-i} (m+N-i)_{N-i} t^k\\
=& \sum_{k=0}^\infty \Big( \sum_{i=0}^{ \left[ \tfrac{N}{2} \right]}  \sum_{m=0}^k {\tfrac{N}{2}-i \choose k-m} (m+N-i)_{N-i}\\
&\times (-4)^{k-m} b_i(N) C_{m+N-i} \Big) t^k.
\end{split}\end{equation}
Thus, from  \eqref{55}, we get the following Theorem.

\begin{thm}\label{thm:4}
For $k=0,1,2\cdots,$ and $N=1,2,3\cdots,$we have 
\begin{equation*}\begin{split}
C_k^{(N+1)} =& \frac{1}{N!} \sum_{i=0}^{ \left[ \tfrac{N}{2} \right]} \sum_{m=0}^k {\tfrac{N}{2}-i \choose k-m} (m+N-i)_{N-i}\\
&\times (-4)^{k-m} b_i(N) C_{m+N-i},
\end{split}\end{equation*}
where $b_i(N)'s$ are as in Theorem \ref{thm:3}.
\end{thm}

\begin{rmk}
Combining \eqref{26} with \eqref{54}, we can show that
\begin{equation}\begin{split}\label{56}
 (1-4t)^{\tfrac{N}{2}}  C^{N+1} =& \sum_{i=0}^{ \left[ \tfrac{N}{2} \right]}  \sum_{j=1}^{N-i} \frac{a_j(N-i)}{N!}b_i(N) (1-4t)^{\tfrac{j}{2}} C^{j+1}\\
 =&\sum_{j=1}^N \sum_{i=0}^{\min (N-j,\left[ \tfrac{N}{2} \right]) }  \frac{a_j(N-i)}{N!}b_i(N) (1-4t)^{\tfrac{j}{2}} C^{j+1}.
\end{split}\end{equation}
\end{rmk}
Equivalently, \eqref{56} can be expressed as
\begin{equation}\begin{split}\label{57}
\sum_{i=0}^{\min (N-j,\left[ \tfrac{N}{2} \right]) }  \frac{a_j(N-i)}{N!}b_i(N) = \delta_{j,N}, \quad (1 \leq j \leq N),
\end{split}\end{equation}
when $\delta_{j,N}$ is the Kronecker delta.
\section{Further remarks}
We start our discussion here with the following expansion of $\sqrt{1+y}$:
\begin{equation}\begin{split}\label{58}
\sqrt{1+y} = \sum_{n=0}^\infty {2n \choose n} \frac{(-1)^{n-1}}{4^n (2n-1)} y^n.
\end{split}\end{equation}
Integrating both sides of \eqref{58} from 0 to 1, we immediately obtain
\begin{equation}\begin{split}\label{59}
\sum_{n=0}^\infty C_n \frac{(-1)^{n-1}}{4^n(2n-1)} = \frac{1}{3}\big(4 \sqrt{2} -2\big).
\end{split}\end{equation}
Next, we integrate the generating function of the Catalan numbers from 0 to $\tfrac{1}{4}$.
\begin{equation}\begin{split}\label{60}
\int_0^\frac{1}{4} \frac{2}{1+\sqrt{1-4t}} dt = \int_0^\frac{1}{4} \sum_{n=0}^\infty \frac{1}{n+1} {2n \choose n} t^n dt
\end{split}\end{equation}
The left hand side of \eqref{60} is, after making the change of variable $t= \tfrac{1}{4} (1-y^2)$, equal to
\begin{equation}\begin{split}\label{61}
\int_0^1 \frac{y}{1+y} dy = \left[ y-\log(1+y) \right]_0^1 = 1-\log 2.
\end{split}\end{equation}
Thus, from \eqref{60} and \eqref{61}, we get the following identity.
\begin{equation}\begin{split}\label{62}
\sum_{n=0}^\infty \frac{1}{(n+1)^2} {2n \choose n} \left( \frac{1}{4} \right)^{n+1} = 1- \log 2.
\end{split}\end{equation}
Finally, again from the generating function of Catan numbers and \eqref{58}, we have
\begin{equation}\begin{split}\label{63}
2 =&\Big( \sum_{l=0}^\infty C_l t^l \Big) \Big( 1+ \sqrt{1-4t} \Big)
\\=&\Big( \sum_{l=0}^\infty C_l t^l \Big) \Big( 1- \sum_{m=0}^\infty {2m \choose m} \frac{1}{2m-1} t^m \Big)
\\=&\sum_{l=0}^\infty C_l t^l - \Big( \sum_{l=0}^\infty C_l t^l \Big) \Big(\sum_{m=0}^\infty {2m \choose m} \frac{1}{2m-1} t^m \Big)
\\=&\sum_{n=0}^\infty C_n t^n - \sum_{n=0}^\infty \left( \sum_{m=0}^n {2m \choose m} \frac{1}{2m-1}C_{n-m} \right) t^n
\\=& \sum_{n=0}^\infty \left( C_n - \sum_{m=0}^n C_m C_{n-m} \frac{m+1}{2m-1}
\right) t^n.
\end{split}\end{equation}
Therefore, from \eqref{63} we obtain the recurrence relation:
\begin{equation}\begin{split}\label{64}
C_n - \sum_{m=0}^n C_m C_{n-m} \frac{m+1}{2m-1} = \left\{ \begin{array}{ll}
2 & \textnormal{if}\,\,\, n=0, \\
0 & \textnormal{if}\,\,\, n>0.
\end{array} \right.
\end{split}\end{equation}
Noting that $C_n = \tfrac{1}{n}{2n \choose n+1}$, we see that \eqref{64} for $n>0$ is equivalent to the following identity.
\begin{equation}\begin{split}\label{65}
{2n \choose n+1} = n \sum_{m=0}^n \frac{m+1}{m(n-m)(2m-1)} {2m \choose m+1} {2n-2m \choose n-m+1}, \quad (n>0).
\end{split}\end{equation}
Further, separating terms corresponding to $m=0$ and $m=n$ from \eqref{64} for $n>0$ and after rearranging the terms, we get the following recurrence relations for the Catalan numbers.
\begin{equation}\begin{split}\label{66}
C_0 =C_1 =1, \quad C_n = \frac{2n-1}{3(n-1)} \sum_{m=1}^{n-1} C_m C_{n-m} \frac{m+1}{2m-1}, \quad (n \geq 2).
\end{split}\end{equation}
Compare \eqref{66} with the recurrence relation in \eqref{02}.

\end{document}